\documentclass[12pt]{article}
\usepackage{amsfonts, amssymb,version}
\usepackage{amsmath}
\def\Gg{{\mathbb G}}

\def\mf{{\mathfrak m}}

\def\P{{\mathbb P}}

\def\Q{\mathbb Q}

\def\Z{\mathbb Z}

\let\lra\longrightarrow

\let\ov\overline
\let\mf\mathfrak

\let\wt\widetilde

\newtheorem{theorem}{Theorem}[section]
\newtheorem{proposition}[theorem]{Proposition}
\newtheorem{lemma}[theorem]{Lemma}
\newtheorem{definition}[theorem]{Definition}
\newtheorem{corollary}[theorem]{Corollary}

\def\rem{\refstepcounter{theorem}\paragraph{Remark \thetheorem}}

\def\proof{\paragraph{Proof}}

\makeatletter \def\l@section{\@dottedtocline{1}{0em}{1.2em}} \makeatother

\begin{document}

\title{Extension of structure groups of principal bundles in positive characteristic}
\author{Fabrizio Coiai and Yogish I. Holla }
\date{}

\maketitle

\begin{abstract}
In this article we study the behaviour of semistable principal
$G$-bundles over a smooth projective variety $X$ under the extension
of structure groups in positive characteristic.  
We extend some results of Ramanan-Ramanathan
\cite{Ramanan-Ramanathan} on rationality of instability flags
and show that the associated vector bundles via representations of
$G$ are not too unstable and the instability can be bounded by a
constant independent of  semistable bundles.
As a consequence of this the boundedness of the set of isomorphism classes
of semistable $G$-bundles with fixed degree and Chern classes is proven.
\bigskip

\it{Mathematics subject  classification 14J60, 14L24}
\end{abstract}


%
\section{Introduction}

Let $G$ be a connected reductive algebraic group over an algebraically
closed field $k$ of arbitrary characteristic. Let $X$ be a smooth
projective variety over $k$ with a  fixed polarization $H$. In this paper we
address the question of what happens to semistability of principal $G$-bundles under the extension of structure groups.

Recall the definition of a rational $G$-bundle $E$ as a principal
$G$-bundle over a big open subscheme (whose complement is of codimension at least $2$). 
A rational $G$-bundle  
$E$ over $X$ is semistable with respect to the polarization $H$ if for
any reduction to a parabolic subgroup $P$ of $G$ over any big open subscheme, 
the line bundle associated to any  dominant character on $P$ has
degree $\le 0$. 

One notes that restrictions of torsion free sheaves to suitable open
sets define rational ${\rm GL}(V)$-bundles and in this case the above
definition of semistability coincides with usual $\mu$-semistability.

Let $\rho :G \lra {\rm GL}(V)$ be a representation of $G$ on a vector
space $V$ which sends the connected component of the center of $G$ to that of ${\rm GL}(V)$.
 For any rational $G$-bundle $E$ we denote by $E(V)$ the associated
 rational vector bundle. 
 
 When the characteristic of the field is zero, it is proved in 
\cite{Ramanan-Ramanathan} that the bundle $E(V)$ is semistable.
If the characteristic of the field is a prime $p$ which is
sufficiently large (quantified by the height of the representation) 
then the semistability of $E(V)$ is proved in
\cite{I-M-P}. 

In the case of arbitrary characteristic it is known that the bundle $E(V)$ need not be semistable.
Let $\mu_{{\rm max}}(E(V))$ (and $\mu_{{\rm min}}(E(V))$) be the slopes of the first (and the last) term
in the Harder-Narasimhan filtration of $E(V)$. We prove the following theorem.

\begin{theorem}\label{main}
Let $\rho:G \to {\rm GL}(V)$ be a representation which sends the 
connected component of the center of $G$ to that of ${\rm GL}(V)$.
Then there exists a constant $C(X\,,\,\rho)$ (depending only on $X$ and $\rho$) such that
for each rational semistable $G$-bundle $E$ over $X$ we have 
$$
\mu_{{\rm max}}(E(V))-\mu_{{\rm min}}(E(V)) \le C(X\,,\,\rho)
$$
\end{theorem}

We briefly describe the proof. Let $E(G)$ be the group scheme
associated to $E$ and $E(G)_0$ be the  group scheme at the generic
point ${\rm Spec}(k(X))\to X$. Let $P$ be a maximal parabolic subgroup of $ {\rm GL}(V)$ 
and let $\sigma$ be a rational reduction of structure group of 
$E({\rm GL}(V))$.
There is an action of $E(G)_0$ on the smooth projective
variety $E({\rm GL}(V)/P)_0$ over $k(X)$ which is linearized by a suitable line
bundle. The section $\sigma$ gives a $k(X)$-valued point $\sigma _0$ of  
$E({\rm GL}(V)/P)_0$. 
It is known that if $\sigma _0$ is a semistable point for the above action
then the reduction $\sigma$ does not violate the semistability of  
$E({\rm GL}(V))$ (see Proposition (\ref{RR1})). Also, if $\sigma_0$ is not semistable and its instability parabolic
$P(\sigma_0)$ (see Section 3 for the definition) is defined over the function field $k(X)$ of $X$
then again $\sigma$ does not violate the semistability of $E({\rm GL}(V))$ (see Proposition (\ref{RR2})).
This argument in characteristic zero proves that $E(V)$ is semistable because $P(\sigma_{0})$ is defined over $\ov{k(X)}$ and
by its uniqueness it is invariant by the Galois group hence it is defined over $k(X)$. 

In the case of characteristic $p$ one of the important points in our proof is to show that 
 there is an integer $N$ (independent of the $G$-bundle $E$)
such that if  $\sigma _0$ is not 
semistable then its instability parabolic is defined over the field
$k(X)^{p^{-N}}$ (see Proposition \ref{independence}). This part is achieved by 
repeated use of an algebraic result which enables us to get uniform bounds for non-reducedness
of fibers of morphisms of algebraic varieties (see Proposition \ref{bound3}).

Once the instability parabolic is defined over $K^{p^{-N}}(X)$, this parabolic gives rise
to a reduction of structure group of the Frobenius pull-back $(F^N)^*E$.
Using this reduction and some geometric invariant theory
arguments, we reduce the problem to proving the following Theorem which bounds instability of Frobenius pull-backs. 
 
\begin{theorem}\label{sun}
There exists a constant $C(X\,,\,G)$ and a constant $N(G)$ such that for
any rational
$G$-bundle $E$  we have 
$$
{\rm Ideg}(F^*E) \ge p\,N(G) \,{\rm Ideg}(E) \,+\, C(X\,,\,G)
$$
\end{theorem}
(The instability degree ${\rm Ideg}$ is defined by equation (\ref{instdeg}) (in Section 2).

In the case of vector bundles, the above result was proved by X. Sun (see
\cite{Sun}) and Shephard-Barron
(see \cite{Shephard-Barron}).

We use the Theorem (\ref{main}) for groups of lower semisimple rank  
to prove the Theorem (\ref{sun}). In fact we prove a generalization of Theorem (\ref{main}) 
where we replace ${\rm GL}(V)$ by an arbitrary reductive group and (\ref{sun}) will then be a special case 
of this result  (see Remark (\ref{general})).

Let ${\wt c}_i \in A^*(X)$ for $i=1 \ldots n$ 
be elements with $n={\rm dim}(X)$.
Let $S_b(r\,;{\wt c}_1\,\dots, {\wt c}_n)$ be the set of isomorphism classes
of torsion free sheaves $V$ of rank $r$ and $c_i(V)={\wt c}_i$
satisfying $ {\mu}_{{\rm max}}(V)\,-\,{\mu}_{{\rm
min}}(V) \,\le\, b $.

Let $c_i \in A^i(X)$ for $2 \le i
\le n$ be fixed. We also fix a homomorphism 
$d \in {\rm Hom}({\cal  X}(G)\,,\,A^1(X))$. Here $A^k(X)$'s are the Chow
groups and ${\cal  X}(G)$ is the group of characters.

In the last section, we use the above results to show
the following result on boundedness of semistable $G$-bundles (which are defined on all of $X$).

\begin{theorem}\label{boundedness}  
Assume that the set $S_b(r\,;{\wt c}_1\,,\ldots \,,{\wt c}_n)$ is bounded
for all choices of ${\wt c}_i$, $b$ and $r$. Then 
the set ${\cal S}_G(d\,;c_2,\ldots,c_n)$ of isomorphism classes of semistable $G$-bundles 
$\{E\}$ with degree $d_E=d$ and $c_i({\rm ad}(E))=c_i$ is bounded.
\end{theorem}  
Here the degree of a principal bundle $E$ is an element 
$d_E\in {\rm Hom}({\cal  X}(G)\,,\,A^1(X))$ 
defined by $d_E(\chi)=c_1(\chi_*(E))$ for any character $\chi$
of the group $G$. Here $\chi_*(E)$ is the line bundle associated to 
$E$ via $\chi$.

In characteristic $0$, the boundedness of the set 
$S_b(r\,;{\wt c}_1,\dots, {\wt c}_n)$ is well known (see \cite{H-L}
for example). In the case of positive characteristic, for surfaces,
this is due to Gieseker \cite{Gieseker} and Maruyama \cite{Maruyama1}. For higher dimensional
varieties this is recently claimed by Langer \cite{Langer}. This along with Theorem (\ref{boundedness}) would
then prove boundedness of semistable $G$-bundles over $X$ with fixed Chern classes.

When $X$ is a smooth projective curve in characteristic $0$, the
boundedness of the semistable $G$-bundles with fixed degree 
is due to Ramanathan \cite{Ramanathan};
in the case of positive characteristic, it is proved in
\cite{Holla-Narasimhan}. \\
\textbf{Acknowledgments} The authors wish to thank M.S. Narasimhan who
initiated them into this problem, and also for sharing his ideas.
The first named author would like to thank the Tata Institute of Fundamental Research
(Mumbai and Bangalore centers) for hospitality.

\medskip

\section{Basic definitions and notations }
\begin{center}
\textbf{Semistable $G$-bundles}
\end{center}
In this section we recall and prove some basic facts about principal $G$-bundles over varieties.
Let $k$ be an algebraically closed field. 
Let $G$ be a connected
reductive algebraic group over $k$. Let $T$ be a maximal torus and $B$ a Borel
subgroup containing $T$. Let $R_u(B)$ be the unipotent radical of $B$.
Then $B$ is a semi-direct product $R_u(B)\cdot T$. 
We denote by ${\cal X}_*(T)$, the group
of 1-parameter subgroups of $T$ (denote by 1-PS). ${\cal X}^*(T)$
denotes the group of characters of $T$.  We have a perfect pairing
${\cal X}_*(T) \otimes {\cal X}^*(T) \lra \Z$ which will be denoted by
$(\cdot \,,\,\cdot)$.  Let $\Phi \subset {\cal X}^*(T)$ be the set of
roots  of $G$, $\Phi ^+$ be the set of positive roots and $\Delta$ be
the set of  simple roots corresponding to $B$. 
For any $\alpha\in \Phi$, let $T_{\alpha}$ be the 
connected component
of ${\rm ker}(\alpha)$ and $Z_{\alpha}$ the centralizer of
$T_{\alpha}$ in $G$. Then the derived group
$[Z_{\alpha},\,Z_{\alpha}]$ is of rank one and there is a unique 1-PS
${\hat \alpha}: \Gg _m \lra T\cap [Z_{\alpha},\,Z_{\alpha}]$ such that
$T=({\rm im}\,{\hat \alpha}) \cdot T_{\alpha}$ and 
$({\hat \alpha},\,\alpha)=2$. This ${\hat \alpha}$ is the coroot 
corresponding to
$\alpha$. We denote by ${\hat \Phi}$ the set of coroots. The quadruple
$\{{\cal X}^*(T),\, \Phi,\,{\cal X}_*(T),\,{\hat \Phi}\}$ defines a
root system. For each $\alpha \in \Delta$ we have the fundamental dominant 
weight $w_{\alpha}\in {\cal X}^*(T)\otimes \Q$ defined by 
$({\hat \beta},w_{\alpha})=\delta _{\alpha, \beta}$ for $\beta \in \Delta$
and $(\gamma,w_{\alpha})=0$ for any 1-parameter group in the connected 
component of the center of $G$.  Let $W= N(T)/T$ be the Weyl group. We fix
a $W$-invariant inner product on ${\cal X}_*(T) \otimes \Q$ (hence on ${\cal X}^*(T) \otimes \Q$).

Let $P$ be a parabolic subgroup of $G$ containing $B$. Let $U$ be its 
unipotent radical.
Then there is a subset $\Pi \subset \Delta$ such that $P=P_{\Pi}$. 
Let $Z_{\Pi}=(\cap _{\alpha \in \Delta-\Pi}{\rm ker} \alpha)^0$ 
be the connected
component of the intersection of kernels of roots in $\Delta
-\Pi$. By taking the centralizer of $Z_{\Pi}$ one obtains a splitting $P \to P/U=L$
with $Z_{\Pi}$ being the connected component of the center of $L$.

Let $X$ be a smooth projective variety over $k$ of dimension $n$ and let $k(X)=K$ be the function field of $X$.
Let $H$ be a fixed polarization on $X$. Since any line bundle $\mathcal{L}$ over a 
big open subscheme $U$ (whose complement is of codimension at least $2$)   
 admits a unique extension to all of $X$, its first Chern class makes sense.
Recall the definition of the degree of the line bundle to be 
${\rm deg}(L)=c_1(L)\cdot H^{n\,-\,1}$. Hence for any torsion free
sheaf its first Chern class and the degree with respect to $H$ makes sense. 

We recall from \cite{Ramanan-Ramanathan} that a \textbf{rational $G$-bundle} $E$ is a principal
$G$-bundle over a big open subscheme of $X$. For $G={\rm GL}(V)$ this defines a
vector bundle over a big open subscheme and we call it rational vector 
bundle. 

Let $E\lra U\subset X$  be a rational $G$-bundle with $U$ a big open subscheme.
By \textbf{rational reduction} of structure group $\sigma$ of $E$ to $P$ we mean a
reduction of structure group over a big open subscheme $U' \subset U$.
More precisely it is a pair
$(E_{\sigma},\phi)$ with $E_{\sigma}$ a $P$-bundle over $U'$ and an 
isomorphism $\phi:E_{\sigma}(G) \lra E|_{U'}$. This is equivalent to giving a section  $\sigma$ of the fiber bundle $\pi: E/P \lra U$ over 
$U'$. Here $E/P$ denotes the extended rational fiber bundle $E(G/P)$ over $X$.

Let $T_{\pi}$ be the tangent bundle along the fibers of the map $\pi$.
Then $T_{\pi}$ is a rational vector bundle. 
For a reduction of structure group $\sigma$ we will denote by $T_{\sigma}$ 
the rational vector bundle defined by the pull-back of $T_{\pi}$ under 
$\sigma$. 
We will also fix notations for the Lie algebras by putting ${\mathfrak
g}$, ${\mathfrak p}$ and ${\mathfrak l}$ for the Lie
algebras of $G$, $P$ and $L$ respectively.
Then it can be verified that $T_{\sigma}$ is the rational vector bundle on $X$ 
associated to $E_{\sigma}$ for the natural representation of $P$ on ${\mf g}/{\mf p}$.

Recall the following definition of semistability from Ramanan-Ramanathan 
\cite{Ramanan-Ramanathan}. 
\begin{definition}
A rational $G$-bundle  
$E \lra U \subset X$, with $U$ a big open set is \textbf{semistable} with
respect
to polarization $H$ if for any reduction of $E$ to any parabolic subgroup 
$P$ of
$G$ over any big open set $U'\subset U$, 
the line bundle associated to any dominant character on $P$ has
degree $\le 0$.
\end{definition}
 This definition is equivalent to the fact that 
${\rm deg}(T_{\sigma})\,>\,0$ for each rational parabolic reduction.

If $V$ is a torsion free sheaf then over a big open set $U$ the
restriction $V|_U$ is
a vector bundle. The above definition of semistability is equivalent
to the $\mu$-semistability of $V$.      
\begin{center}
\textbf{Instability degree and Harder-Narasimhan reduction}
\end{center}
If the rational $G$-bundle is not semistable then
there is a notion of Harder-Narasimhan reduction which we recall here.
For a rational $G$-bundle $E$ which is not semistable we define the \textbf{instability degree} to be:
\begin{equation}\label{instdeg}
{\rm Ideg}(E) = {\rm Min}_{\{P,\sigma \}}{\rm deg}(T_{\sigma})
\end{equation}
where the minimum is taken over all parabolic subgroups $P$ and rational
reductions $\sigma$. If the rational $G$-bundle is semistable then we
say its instability degree is $0$.
 The following lemma shows that the instability degree makes sense, and is an analogue of Lemma 2.1 of 
\cite{Holla-Narasimhan} for the higher dimensional varieties.

\begin{lemma}
There exists a constant $A_E$ such that for any rational reduction
$\sigma$ of $E$ to any parabolic $P$ we have ${\rm deg}(T_{\sigma})\,>\,A_E $
\end{lemma}
\proof
It is enough to show that the degree of the rational vector subbundle 
${\rm ad}(E_{\sigma}) \subset {\rm ad} (E)$ is bounded above. 

We can first extend the bundle ${\rm ad} (E)$ to get a torsion free
sheaf ${\cal E}$ on $X$. Then we can extend ${\rm ad}(E_{\sigma})$ inside
${\cal E}$ to obtain a torsion free subsheaf.
There exists a constant $A_E'$ such that for  
any curve $C$ in the class $|H^{n-1}|$, we have a bound
$h^0(C\,,{\cal E}|_C) \,\le \,A_E'$. Let $g$ be the maximum of the
genus of smooth curves in $|H^{n-1}|$.
 Now if $C$ is a smooth
projective curve which sits in the domain of definition of 
${\rm ad}(E_{\sigma})$ and ${\rm ad} (E)$ then we get 
${\rm deg}(E_{\sigma}) \,\le\,A_E'\,+\, (g-1)\, {\rm rank}({\rm
  ad}(E))$.
This proves the lemma.  $\hfill \square$

\begin{definition}
A rational reduction of structure group $\sigma$  of $E$ to a parabolic $P$ is said to be a
\textbf{Harder-Narasimhan reduction} if 
${\rm deg}(T_{\sigma})={\rm Ideg}(E)$ and $P$ is maximal among
parabolic subgroups of $G$ containing $B$ for which the above equality
holds.
\end{definition}

The Harder-Narasimhan reductions as defined above satisfy the
following properties stated in Ramanathan \cite{Ramanathan78}. (see \cite{Biswas-Holla} for a proof).

\begin{enumerate}
\item If $L$ is the Levi quotient of $P$, then the principal
$L$-bundle $E_{\sigma}(L)$ obtained by extending the structure
group is semistable;

\item After fixing a Borel subgroup $B\, \subset\, P$ of $G$,
for any nontrivial character $\chi$ of $P$ which is a
non-negative linear combination of simple roots, the
associated rational line bundle $\chi _*((E_{\sigma}))$ over $X$
is of positive degree.
\end{enumerate}

It is proved in Behrend \cite{Behrend} that over a smooth projective
curve there is a unique reduction to a parabolic subgroup 
containing $B$ satisfying the above properties. In the case when $X$ is higher dimensional the uniqueness is
known only when the characteristic of the field is $0$ or it is a
large prime $p$.

We will not have the occasion to use the uniqueness of the above 
reduction. We will only use its existence.

For the case $G={\rm GL}(V)$ the above reduction defines the
Harder-Narasimhan filtration and the uniqueness is then immediate. We have the following lemma
which compares the instability degree with the 
$\mu _{{\rm max}}\,-\,\mu _{{\rm min}}$ of the rational vector bundle.

\begin{lemma}\label{GLN}
Let $E$ be a rational principal ${\rm GL}(V)$ bundle of rank $r$
over $X$ which is not semistable (we will denote by $E(V)$, the associated
vector bundle). Then we have the following.
$$
\mu _{{\rm max}}(E(V))\,-\,\mu _{{\rm min}}(E(V)) \,\le \,
-\frac{2}{r^2} \, {\rm Ideg}(E)
$$
\end{lemma}
\proof
For the proof one first notices that if $F \subset E(V)$ is a rational
subbundle of rank $r_1$ and $F_1$ is the quotient, then it defines a 
rational reduction of
structure group $\sigma$ of $E$ to a maximal parabolic $P_1$ of ${\rm GL}(V)$. One further has
an isomorphism of rational bundles $T_{\sigma}\cong {\cal H}om(F\,,\,F_1)$.

This implies that 
$ \mu(F)-\mu(F_1)=- \mu ({\cal H}om(F\,,\,F_1))\,\le 
\,-{\rm Ideg}(E)/(\,r_1\,(r-r_1))$.
 
This inequality can also be written by eliminating $F_1$ or $F$,
we get $ \mu(F)-\mu(E(V)) \,\le\,-{\rm Ideg}(E)/\,( r\,r_1\,)$ and 
$ \mu(E(V))-\mu(F_1) \,\le\, -{\rm Ideg}(E)/\,(r\,(r\,-\,r_1))$.

Now if we take $F$ to be the rational subbundle which is maximal
destabilizing then we have $\mu _{{\rm max}}(E(V))=\mu(F)$. This implies that 
$\mu _{{\rm max}}(E(V)) \,-\,\mu(E(V)) \,\le\, -{\rm Ideg}(E)/\,r^2$. Similarly
one has  $\mu(E(V))\,-\,\mu _{{\rm min}}(E(V))  \,\le\, -{\rm
Ideg}(E)/\,r^2$. Combining these we have the proof of the lemma.
$\hfill \square$
\begin{center}
\textbf{Frobenius morphism}
\end{center}
Let $k$ be an algebraically closed field of characteristic $p>0$.  
The $p^m$-th power map $\mathcal{O}_X \to \mathcal{O}_X$ 
given by $f \to f^{p^m}$ is a homomorphism and gives rise to a morphism 
$F^m_X:X\to X$ called the absolute Frobenius morphism. 
Using this morphism we can pull-back a $G$-bundle $E$ over $X$ to get a 
$G$ bundle $(F_X^m)^*(E)$ which we will denote by $(F^m)^*E$. 
One notes that the 
absolute Frobenius morphism is not a $k$-morphism. But since $k$ is a perfect
field we can always twist by the Frobenius isomorphism of $k$ to ensure that the $G$-bundle
$(F^m)^*(E)$ has a well defined  $k$-structure. 

We will use the following well known result about Frobenius morphisms. 

\begin{proposition}\label{connection}
Let $E$ be a rational $G$-bundle over $X$. Then there exists a $p$-connection 
$\nabla$ on the 
$G$-bundle $F^*(E)$ which satisfies the following property: 
for any rational reduction of 
structure group $\sigma$ of  $F^*(E)$ to a parabolic $P$ there is a  
vector bundle map (second fundamental form) 
$\nabla _{\sigma}:T_X \to T_{\sigma}$ (wherever $T _{\sigma}$ is defined) 
such that 
the following are equivalent.
\begin{enumerate}
\item There exists a rational reduction $\sigma _0$ of $E$ to $P$ such that  
$F^*(\sigma _0)=\sigma$. \\ 
\item $\nabla _{\sigma}$ is zero.
\end{enumerate}
\end{proposition}
\proof See for example, proof of Theorem (2.1) of \cite{Sun}.

\section{Geometric invariant theory and the method of Ramanan and Ramanathan}

In this section we will describe some basic facts about geometric
 invariant theory and briefly explain the results of Ramanan and Ramanathan 
 \cite{Ramanan-Ramanathan}.

\begin{center}
\textbf{The instability parabolic}
\end{center}

Let $K$ be a field and let $ {\ov K}
$ be a fixed algebraic closure.
Let $G$ be a connected reductive algebraic group over $K$. Let $V$ be a finite dimensional
 representation of $G$ we get an induced action of $G$ on the projective space $\P(V)$ of lines in $V$.
For a point $v \in \P(V)$ we will denote, by abuse of notation, a representative in $V$ again by $v$.

Firstly we will describe the theory when $K={\ov K}$ is 
algebraically closed and later extend the theory to non-algebraically closed fields.

	Recall that a point $0 \neq v \in V$ is \textbf{semistable} for the $G$-action if the
	closure $ { \ov{G \, v}}$ of the orbit of $v$ does not contain $0$.
One knows that this definition is equivalent to existence of a 
$G$-invariant element $\phi \in S^n(V^*)$ for some $n> 0$ such that 
$\phi(v)\neq 0$ .

If $v\in V$ is not semistable then recall the following notions.

For a 1-PS $\lambda$ of $G$, consider the decomposition of $V= \bigoplus V_i$ with $V_i=\{ v\in V\,|\,\lambda(t)v=t^i\,v\}$.
For an $v\in V$ one defines the invariant 
$$
m(v\,,\lambda)\,=\,{\rm inf}\{ i\,|v \,\,\mbox{has a non-zero component in} 
\,\,V_i\}
$$

 Using the $W$-invariant inner product on a fixed maximal
    torus $T$ one defines the slope 
    $ \nu (v \,, \lambda)\, = \, m (v \,, \lambda) / ||\lambda || $ 
    for all 1-PS in the maximal torus $T$. Since maximal tori are conjugates
    this definition can be extended to all 
    $1$-PS in $G$.

The following lemma will be used in the sequel.
\begin{lemma}\label{git1}
With the above notations there exists a constant $C$ (independent of $v\in V$ and 
$1$-PS $\lambda$) such that $\nu(v\,,\lambda) \le C$.
\end{lemma}
\proof   See Proposition (2.17) p. 64 of \cite{Mumford} for a proof. 

We define the \textbf{instability $1$-PS} for a given 
$v \in V$ (which is not semistable) as one for which $ \nu(v\,,\lambda)$ attains its maximum among all $1$-PS of $G$ (see Theorem (1.5,\,a) of
\cite{Ramanan-Ramanathan}). 

For a $1$-PS $\lambda$, recall the definition of the parabolic 
$P(\lambda)$ whose valued points are characterized by elements $g \in G$ for which the limit
${\rm Lim}_{t \lra 0}\lambda(t)\,g\,\lambda(t)^{-1}$ exists.

In the following Proposition we summarize the basic facts in 
geometric invariant theory.

\begin{proposition}\label{git2} 
Let $v \in V$ be a non-zero element which is non-semistable.
\begin{enumerate}
\item There is a unique parabolic subgroup 
$P(v)$ with the property that for any instability $1$-PS $\lambda$ for $v$ 
we have $P(v)= P(\lambda)$. \\

\item For any maximal torus $T \subset P(v)$ there is a unique $1$-PS 
$\lambda _T \subset T$ such that it is an instability $1$-PS for $v$.
\end{enumerate}
\end{proposition}
\proof
See Theorem (1.5,\,b and c) \cite{Ramanan-Ramanathan}.

The above uniquely defined parabolic $P(v)$ will be called the \textbf{instability 
parabolic} for $v$. Here the uniqueness of $\lambda$ is as a subgroup of $T$ rather than a morphism $\Gg_m \to T$.

If $G$ acts on projective variety $M$ which is linearized  by an ample line bundle $\cal{L}$
then by taking some power of $\cal{L}$ we get a representation $V$ of $G$ and a $G$-equivariant 
embedding $i:M \to \P(V)$ with $i^*{\cal{O}}(1)$ being some power of $\cal{L}$. In this 
setup we say a point $m \in M$ is semistable for $G$ action if the corresponding point in $V$ is semistable.

Let $v \neq 0 $ be a non-semistable point. Let 
$P \, = \, P (v)$ be its instability
    parabolic and let $ \lambda \subset T \subset P$ be a chosen tuple
    of instability $1$-PS and a maximal torus $T$. Let 
    $V \, = \, \bigoplus _i V_i$ be the decomposition of $V$ with 
    respect to $ \lambda$. Let $j \, =\, m (v \,, \lambda) $. 
    Using this we have a decomposition $v \, = \sum_{i \ge 0} v_i$, 
    where $v_i \in V_ {i + j}$.

Here one notes that $V ^ j = \bigoplus _ {i \ge j} V_i$ is preserved 
under the action of $P=P (v)$ and the unipotent radical $U \subset P$ 
pushes $V ^ j$ to $V ^ {j + 1}$, thus giving an action
    of the Levi quotient $L \, = \, P / U$ on $V ^ j / V ^{j + 1}$.

The $W$-invariant inner product on a fixed maximal torus of $G$ naturally gives rise to a 
$W$-invariant inner product
 on $T\subset P(v)$. Let $l_{\lambda} \in {\cal X}^*(T)\otimes \Q$ 
 be the dual of $\lambda$. Let $r_1 \in \Z^{+}$ such that $r_1\,l_{\lambda}$ 
 defines a character of $T$. The restriction of this character to the connected component $Z^0(L)$of the center of $L$, and taking a further multiple, extends to give a character of $L$.
Hence given a $\lambda $ we get a character 
 $\chi$ of $P(v)$ which is well defined up to a positive integral multiple.
 
In the following proposition we describe the basic result of Ramanan-Ramanathan \cite{Ramanan-Ramanathan}
concerning the behaviour of $v_0\in V^j/V^{j+1}$ under
the induced action of $L$ on $V^j/V^{j+1}$.  
\begin{proposition}\label{git3} 
Assume that the group $Z^0(G)$ acts trivially on $V$. Then there exists a positive 
integer $r$ and dominant character $\chi$ of $P$ such that the point 
$v_0 \in \P(V^j/V^{j+1})$ is semistable for the natural action of $L$ 
with respect to 
linearization given by ${\cal O}(r)\otimes {\cal O}_{\chi ^{-1}}$, where 
${\cal O}_{\chi ^{-1}}$ is the trivial line bundle with $L$ acting on it 
by $\chi ^{-1}$
\end{proposition} 

The proof of the above result (as given in Proposition (1.12) in \cite{Ramanan-Ramanathan}) also gives a recipe 
to find the integer $r$ and the character $\chi$ and they are related by the following Lemma.

\begin{lemma}\label{git4}
There is a character $\chi '$ of the maximal torus $T \subset P$ 
such that the following holds.
\begin{enumerate}
\item $\chi'|_{Z^0(L)}\,=\,\chi |_{Z^0(L)}$
\item $\chi' \,=\,r\,\nu(v\,,\lambda)\,||\lambda || \,l_{\lambda}$, where
$l_{\lambda}$ is the dual of $\lambda$.
\end{enumerate}
 \end{lemma}

The above Lemma will be used in our proof of the main Theorem.

\begin{center}
\textbf{The Rationality of the instability parabolic}
\end{center}

We will now assume that the ground field $K$ is not algebraically closed.  
Let $G$ be a connected reductive group over $K$ which
acts on a projective $K$-variety $M$, linearized by an ample line bundle $\cal{L}$, 
thus giving a $G$-equivariant embedding $i:M \to \P(V)$ as before.

We will call a $K$-valued point $v \in V$ semistable if it is so after a 
base change to algebraic closure. In this way we will avoid the confusion of 
which field the semistability definition is used.

Let $m$ be a $K$-rational point of $M$ which is not semistable. Let $P(m)$ 
be its instability parabolic defined over ${\ov K}$. 

\rem \label{uniqueness} Note that if $P(m)$ is defined over
$K_s$ then it is already defined over $K$. This is 
because of the uniqueness of $P(m)$ (see Proposition \ref{git2} (1)) 
and the Galois descent argument.
Also note that if $P(m)$ is defined over $K$ then it contains a maximal torus 
over $K$ which splits over $K_s$. Then the instability $1$-PS of $m$
which is contained in the maximal torus over $K_s$, by uniqueness 
(Proposition \ref{git2} (2)) is Galois invariant and hence it is defined over $K$.

Let $O(m)$ be the (reduced) orbit of $G$
at $m$. Since $m$ is defined over $K$ the orbit $O(m)$ is also defined over $K$.

We briefly recall the construction of a scheme $M(P)_{x_m}$ which will be 
used later in an important way.

We can find a $g \in G$ such that $g\,P(m)\,g^{-1}=P$ is defined over
$K_s$ (as the variety of all parabolics conjugate to $P(m)$ is
defined over $K_s$, being absolutely reduced,  has a $K_s$-rational point
(see \cite{Borel-Tits})). If $x_m=g\,m$ then $P=P(x_m)$ is  the
instability parabolic of $x_m$. 

Since $P$ is defined over $K_s$ and over this field $G$ splits, we have a 
maximal torus in $P$ which splits over $K_s$. 
Hence there is a instability $1$-PS $\lambda$ of $x_m$ in this maximal 
torus over $K_s$.   

The representation $V$ of $G$ decomposes as $V= \bigoplus _{i\in \Z}V_i$ for the action
of $\lambda$, where 
$V_i\, = \, \{v \in V|\lambda(t) \cdot v \,= \, t^i \,v, t\in \Gg_m\}$.
Let 
$j\,= m(x_m,\, \lambda)$
and $V^j =\bigoplus_{i \ge j}V_i$. 

Recall the definition of the $K_s$-scheme $M(P)_{x_m}$ as
the scheme theoretic intersection of the $K_s$-subschemes 
$\P (V^j)$ and $O(m)$ of $\P(V)$. 

The following two results summarizes the basic properties of the scheme $M(P)_{x_m}$.

\begin{lemma}\label{gitlemma}
The ${\ov K}$-rational points of the $K_s$-subscheme $M(P)_{x_m}$ 
of the $K$-scheme $O(m)$ are precisely those points which have 
$P(x_m)$ as their instability parabolic. Moreover, when the $G$
action on $m$ is strongly separable then $M(P)_{x_m}$ is absolutely
reduced. 
\end{lemma}
\proof
See Lemma (2.4) of \cite{Ramanan-Ramanathan}.

Recall that the $G$-action at $m \in M(\ov{K})$ is said to be strongly separable 
if the isotropy subgroup scheme $G_x$ is reduced at every point $x\in M({\ov K})$ which is 
in the closure of the orbit $O(m)$.

\begin{lemma}\label{parabolic}
Suppose that $y\in M(P)_{x_m}\subset O(m)$ is a $K_s$-rational point and that there is an $h\in G(K_s)$ 
such that  $h$ maps to $y$ under the orbit map $G \to O(m)$. Then $P(m)$ is defined over $K$.
\end{lemma}
\proof   
By lemma (\ref{gitlemma}), the point $y =h\, m$ has the property that $P(y)=P(x_m)$, hence $P(y)$ is defined over $K_s$.   
This implies that $P(m)=h\,P(y)\,h^{-1}$ is also defined over $K_s$, hence by Remark (\ref{uniqueness}) 
 we conclude the  proof of the Lemma. $\hfill \square$

The above lemma has the consequence that if the action of $G$ is strongly
separable at $m$ then the parabolic $P(m)$ is already defined over
$K$ (also see Proposition 2.4, \cite{Ramanan-Ramanathan}).

\begin{center}
\textbf{The argument of Ramanan and Ramanathan }
\end{center}

Let $E$ be a rational $G$-bundle over $X$.
Let $\rho: G \lra G_1$ be a representation of $G$ which takes the connected component of the center of $G$ to the center of $G_1$. 
Let $P_1$ be a parabolic
subgroup of $G_1$.  We fix a
representation  $ G_1 \lra {\rm GL}(V_{P_1})$ such that it
defines an embedding of $ G_1/P_1 \subset \P(V_{P_1})$ with the
property that the character of $P_1$ on $V_{P_1}$ is a positive 
multiple $m_{P_1}$ of the character of $\chi_{P_1}$ associated to the
restriction of the adjoint representation of $P_1$ on the vector space 
${\mf g}_1/{\mf p}_1$.

The line bundle ${\cal O}(1)$ on the projective variety $\P(V)$ gives
rise to a line bundle $\cal{L}$ over $G_1/P_1$.

For the rational $G$-bundle $E$ over $X$ we have the group
scheme $E(G)$ using the conjugation action of $G$ on itself. Let $E({\cal{L}})$ over
$E(G_1/P_1)$ be the associated line bundle over the associated 
rational fiber bundle over $X$.

 Let
$E(G)_0$ be the group scheme defined over the function field $K$ of $X$.
We also have the action of $E(G)_0$ on the projective variety 
$E(G_1/P_1)_0 \subset E(\P(V))_0$ over $K$, linearized by the line bundle $E({\cal{L}})_0$.

 Let $\sigma$ be a
reduction of structure group of $E(G_1)$ to $P_1$. Let $\sigma _0$
be the associated $K$-rational point of $E(G_1/P_1)_0$.

In the following two propositions we summarise the basic argument of Ramanan-Ramanathan.
\begin{proposition}\label{RR1}
Let $\sigma _0$ be semistable for the action of $E(G)_0$ on 
$E(G_1/P_1)_0$ (over ${\ov K}$) for the polarization $E({\cal{L}})_0$. 
Then the section $\sigma$ has the property that ${\rm deg}(T_{\sigma}) \ge 0$.
\end{proposition}
\proof
See Proposition (3.10,\,(1)) of \cite{Ramanan-Ramanathan}.

Suppose that $\sigma _0$ is not semistable for the above action and that the 
instability parabolic $P'$ for $\sigma _0$ is defined over the field $K$. Then
the parabolic $P' \subset E(G)_0$ gives rise to a rational reduction of the structure group
$\tau$ of $E$ to the parabolic $P$ such that $(E_{\tau}(P))_0=P'$.

The following result is slightly more general than the Proposition (3.13) of \cite{Ramanan-Ramanathan} (without the semistability assumption on $E$)  and its
proof is along the same lines.
\begin{proposition}\label{RR2}
Suppose $\sigma _0$ is not semistable and its instability parabolic is 
defined over the field $K$ then there exists 
a positive integer $r$ and a dominant character $\chi$ of $P$ 
(related by the  Lemma \ref{git4}) such that the following inequality holds
\begin{equation}
-(r\, m_P)\,{\rm deg}(T_{\sigma}) \le {\rm deg}(\chi _*(E_{\tau}))
\end{equation}
\end{proposition}

The above result when $E$ is semistable implies that ${\rm{deg}}(T_{\sigma})\geq0$ and this
along with Proposition \ref{RR1} is used in characteristic zero to show
that $E(G_1)$ is semistable.

We will use Proposition \ref{RR2} in this generality because the instability parabolic
will be defined after a suitable Frobenius pull-back of $E$ which may not be semistable (see proof of Theorem \ref{main}).

\section{A result on instability parabolics}

One of the main steps in our proof of Theorem (\ref{main}) is a result (Proposition (\ref{independence}))
which gives uniform bounds for the domain of definition of instability parabolics.
For proving this result we need to estimate the non-reducedness  
of the fibers of morphisms of algebraic varieties and we will do this part  
first.

We start with some definitions which will be used later.
Let $K$ be a field and let $\ov{K}$ be its algebraic closure. 

We define the \textbf{radical index} ${\rm Ri}(A)$ of an affine algebra $A$
over $K$ to be the smallest integer $n$ such that for any 
$f$ in the radical ${\rm Rad}(A)$ of $A$ we have $f^n=0$. For an affine morphism 
$f: Y \lra X$ of finite type $\ov{K}$-schemes we define the radical index 
${\rm Ri}(x)$ of a point $x \in X$ to be ${\rm Ri}(Y_x)$ where $Y_x$
is the fiber of $f$ at the point $x \in X$. 

\begin{proposition} \label{bound2}
Let $f:Y \lra X$ be a morphism of finite
type affine schemes over $\ov{K}$. 
There exists an integer $n$ such that 
${\rm Ri}(x) \le n$ for each $x\in X$.
\end{proposition}

\proof
The proof of this proposition is a series of reductions from the case
of arbitrary $X$ and $Y$ to very specific ones using the induction on
the dimension of $X$. Let $X={\rm Spec}(A)$ and $Y={\rm Spec}(B)$ and $i$
be the homomorphism $A \lra B$. For any prime ideal ${\mf p} \in {\rm Spec}(A)$ (or ${\rm Spec}(B)$)
we write ${\rm Ri}({\mf p}\,,\, B)$ for the radical index of $B/{\mf p}B$.
 
First we may assume that $A$ is integral. This part is an elementary check. 

We use induction on the
${\rm dim}(A)$. So, successively we reduce to the situations where we need to bound
${\rm Ri}({\mf m}\,,\, B)$ for maximal ideals ${\mf m} \in {\rm Spec}(A_f)$ for suitable choices of $f \in A$.
 
We make some reductions on $B$. We may assume that $B$ is
reduced. For this if ${\mf m}$ is a
maximal ideal of $A$ then we can check that
$$
{\rm Ri}({\mf m}\,,\,B) \le {\rm Ri}({\mf m}\,,\,B_{\rm red})\,
+\,{\rm Ri}(B) 
$$

Next we may assume that $B$ is irreducible. Let ${\mf p}_i$, for 
$i=1\ldots m$, be the set of minimal prime ideals in $B$. Let 
${\mf m}$ be a maximal ideal in $A$. We will show that 
$$
{\rm Ri}({\mf m}\,,\,B)\le m \, 
{\rm Max}_{i=1}^m {\rm Ri}({\mf m}\,,\,B/{\mf p}_i)
$$
For this one observes that if $x \in {\rm Rad}(B/{\mf m}B)$ then   
the image of $x$ in each of $B/{\mf p}_i\otimes A/{\mf m}$ lies in 
${\rm Rad}((B/{\mf p}_i)/{\mf m}(B/{\mf p}_i))$. Hence if 
$n={\rm Max}_{i=1}^m {\rm Ri}({\mf m}\,,B/{\mf p})$ then $x$ has
the property that $x^n \in \cap _{i=1}^m({\mf p}_i+{\mf m}B)$. We
can write   $x^n= y_i +z_i$ such that $y_i \in {\mf p}_i$ and 
$z_i \in {\mf m}B$.  Since $B$ is reduced we have   
$\Pi_{i=1}^m (x^n-z_i) \, =\, \Pi y_i =0$, and so 
$x^{nm} \in {\mf m}B$. This proves the assertion. 

Hence from now on we may assume that $A$ and $B$ are integral domains.

We reduce this problem to an open subscheme of $Y={\rm Spec}(B)$. 
Let $b\in B$.
Then there exists an element  $a\in A$ such that $(B/bB)_a$ is flat over  
$A_a$. Hence for any maximal ideal ${\mf m}$ in $A_a$ 
we have ${\rm Tor}^A_1((B/bB)_a\,,\, A/{\mf m})=0$. 

Consider the map 
$B_a \lra B_{ab}$. We will show that 
${\rm Ri}({\mf m}\,,\,B_a) \le {\rm Ri}({\mf m}\,,\,B_{ab})$.
If $x \in B$ is such that some power of it lies in 
$B_a/{\mf m}B_a$ then the image of $x$ also has the same property over 
$B_{ab}$. Then by clearing denominators if 
$r={\rm Ri}({\mf m}\,,\,B_{ab})$ then there is an $m$ such that 
$b^m\,x^r \in {\mf m}B_a$. Using the exact sequence 
$$
\cdots \lra {\rm Tor}^A_1((B/bB)_a\,,\, A/{\mf m})\lra B_a/{\mf m}B_a \lra 
B_a/{\mf m}B_a 
$$
with the last map being multiplication by $b$, we conclude, by the
vanishing of the ${\rm Tor}^A_1((B/bB)_a\,,\, A/{\mf m})$, that $x^r
\in {\mf m}B_a$. 
Hence it is enough to bound radical
index of fibers of some open set of the type $B_b$ over $A$.

We use Noether Normalization (and inverting an element of $A$) to get  
an inclusion $A \hookrightarrow A[x_1\,,\ldots \, ,x_r]=A'
\hookrightarrow B$ such that $B$ is finite over $A'$. 
Let $K_A$(respectively $K_B$ and $K_{A'}$) be the function fields of 
$A$ (respectively $B$ and $A'$). 

Let $L$ be the separable closure of 
$K_{A'}$ in $K_B$. The extension $L \subset K_B$ is purely
inseparable. Hence there is an integer $n_1$ such that for any $x\in
K_B$, we have $x^{p^{n_1}} \in L$. 

Let $C=B\cap L$. Then we observe that $B$ is integral
over $C$ and is a finitely generated $A$-algebra. This implies that $C$ is also a finitely generated $A$-algebra. 
Again by localizing $A$ at an element we
can assume that the $A$ module $B/C$ is flat over $A$ and hence for
any ${\mf m}$ in $A$ we have ${\rm Tor}^A_1(B/C\,,\, A/{\mf m})=0$.
This has the effect that for each maximal ideal ${\mf m}$ in $A$ we
have an injection $C/{\mf m}C \lra B/{\mf m}B$.
Using this and the fact that  $x^{p^{n_1}} \in C$ for
any $x\in B$ we have
$$
{\rm Ri}({\mf m}\,,\,B) \le {\rm Ri}({\mf m}\,,\,C)\,+\,p^{n_1}
$$

The problem now reduces to proving the proposition for the case when 
$A \hookrightarrow B$ is an extension of finitely generated domains
such that the function field extension is separable. Further, it is
enough to prove the result for the case $A \hookrightarrow B_b$ for
some $b\in B$.

We may now assume that $A$ and $B$ are smooth domains. Hence we
conclude that there exist elements $b\in B$ and $a\in A$ such that the morphism 
${\rm Spec}(B_{ab}) \lra  {\rm Spec}(A_a)$ is smooth. This implies that 
the fibers here are reduced and hence proof of the proposition is 
complete. $\hfill \square$

Let $T$ be a finite type scheme over $\ov{K}$. Let 
${\cal N}$ be its radical ideal sheaf. This has the property that 
for any closed point of $T$, the stalk of  ${\cal N}$ is the radical of the
local ring. We define the radical index $\rm{Ri}(T)$ of $T$ to be the smallest integer $n$ such that  
 for any open subset $U \subset T$ and for any
$g\in \Gamma(U\,,\,{\cal N})$ we have $g^n =0\in \Gamma(U\,,\,{\cal N}^n)$.

Let $f:Y \lra X$ be a morphism of finite type $\ov{K}$-schemes. In this
general setting we define the radical index of the closed point $x\in X$
by ${\rm Ri}(x)={\rm Ri}(Y_x)$, where $Y_x$ is the fiber at $x$.
In this case we have the following result which generalizes
Proposition (\ref{bound2}) and this will also be used in our proof of main results.

\begin{proposition}\label{bound3}
Let $f:Y \lra X$ be a morphism of finite type schemes over $\ov{K}$
There exists an integer $n$ such that 
${\rm Ri}(x)\le n$ for all closed points $x\in X$
\end{proposition}
\proof
    Let $\{U_i \}_{i=1}^r$ be an open cover of $X$ by affine open subschemes
with $U_i ={\rm Spec}(A_i)$. Then it is enough to prove the result for the
 case of each of $U_i$,
  hence we may assume that $X={\rm Spec}(A)$ is an affine scheme.

Let $\{V_i\}_{i=1}^{r'}$ be an open cover of $Y$ by a finite number of affine
open subschemes. We write $V_i={\rm Spec}(B_i)$. By Proposition 
(\ref{bound2}), we have positive integers $n_i$ such that 
for each maximal ideal ${\mf m}$ of $A$ we have 
${\rm Ri}({\mf m}\,,\,B_i)\le n_i$. Let $n={\rm Max} \{n_i\}$. 

Then we have ${\rm Ri}(x,Y) \le n$ for each closed point $x\in
X$. This is because the fiber $Y_x$ can be covered by affine open
subschemes $\{V_{i,x}\}$. Here  $V_{i,x}$ is the fiber of $x$ in $V_i$.
 If $U \lra Y_x$ is any open immersion and 
$\nu \in \Gamma(U\,,\,{\cal N})$ then restriction $\nu _i$ of $\nu$ to 
$U\cap V_{i,x}$  lies in $\Gamma(U\cap V_{i,x}\,,\,{\cal N})$. Hence
  the result would follow if we show that for any $K$-algebra $B$ and
  an element $b\in B$ we have  $ {\rm Ri}(B_b) \le {\rm Ri}(B)$.
The last statement is a straight forward verification.
This completes the proof of Proposition (\ref{bound3}). $\hfill \square$ 

\rem One notes that the constant $n$ as defined in the above proposition 
depends on $X$, $Y$, and $f$ but not on the $\ov{K}$-valued points of $X$.

Let $G$ be a reductive algebraic group acting on a variety $M$ (over $\ov{K}$).
For any $x \in M(\ov{K})$ we denote the isotropy subgroup scheme at $x$ by $G_x$. 
The following result is a consequence of the above proposition.

\begin{proposition}\label{bound}
There exists an $N_1$ such that ${\rm Ri}(G_x) \le N_1$ for each $x\in M(\ov{K})$. 
\end{proposition}
\proof
Consider the map $G\times M \lra M\times M$ defined by $(\rho\,, pr_2)$  
where $\rho$ is the action map and $pr_2$ is the second
projection map.
Let $\Delta _M$ be the diagonal map $M\lra M\times M$. Let 
$H = (G\times M)\times _{M \times M} M$. Then we have a natural
projection map $\pi :H \lra M$ which has the property that for any
$x\in M(\ov{k})$ the fiber of the map $\pi$ at $x$ is the isotropy subscheme
$G_x$. The result follows from Proposition (\ref{bound3}). $\hfill \square$

Let $K$ be an arbitrary field and $K_s$ and ${\ov K}$ be its separable
closure and the algebraic closure respectively (in
fact $K$ will be the function field of the smooth projective variety $X$). 

In this case the radical index of a finite type scheme $T$ over $K$ is
defined to be the radical index of the scheme 
${\ov T}=T\otimes _K {\ov K}$.

Let $G$ be a reductive group over $K$. Let $M$ and $V$ be as defined before (in Section 3).

Let $m$ be a  non-semistable $K$ valued point of $M$. 
Let $P(m)$ be the instability parabolic defined over ${\ov K}$.
Recall from Remark (\ref{uniqueness}) that if $P(m)$ is defined over
$K_s$ then it is already defined over $K$. 
Hence $P(m)$ is
always defined over a finite purely inseparable extension of $K$. 

The following Proposition is the main result of this Section.

\begin{proposition}\label{independence} There exists an integer $N$ 
such that for any $K$-rational point $m$ of $M$
which is not semistable, the instability flag $P(m)$ is defined over 
$K^{p^{-N}}$.
\end{proposition}
\proof 

It is enough  to show that there exists an $N$ such that 
the  instability parabolic for any non-semistable $K_s$ rational point of
$M$ is defined over $K_s ^{p^{-N}}$ (see Remark (\ref{uniqueness}).
This enables us to assume that all our objects are defined over the field $K_s$.

Let $m$ be a $K_s$-valued point of $M$ which is not semistable
for the action of $G$ on $M$. Let $O(m)$ be the (reduced) orbit of $G$
at $m$. Let $P(m)$ be its instability
parabolic over ${\ov K}$.
We can find a $g \in G$ such that $g\,P(m)\,g^{-1}=P$ is defined over
$K_s$ (this was seen before). If $x_m=g.m$ then $P=P(x_m)$ is  the
instability flag of $x_m$. 

To this setup we have the scheme $M(P)_{x_m}$ as defined in Section 3 
satisfying the properties of the Lemmas (\ref{gitlemma}) and ({\ref{parabolic}).

Our main goal is to estimate the non-reducedness of the scheme $M(P)_{x_m}$ 
independent of $m$, and this will enable us to prove that 
the parabolic $P(m)$ is defined over a fixed purely
inseparable extension of $K_s$ for all $m$.

Note that we have made a choice of $x_m$ and this choice fixes the
instability $1$--PS $\lambda$ of $x_m$ which is defined over $K_s$. 
Hence the point $x_m$ determines the vector subspace $V^j \subset V$ 
(as defined in Section 3) .

We will show that there exists a positive integer $N_2$  such that 
for any non-semistable point $m \in M$ and any choice of $x_m$ as above, the radical index 
${\rm Ri}(M(P)_{x_m}) \le N_2$. The basic idea of the proof is to prove
that the spaces $M(P)_{x_m}$  occurs as suitable subschemes of the 
fibers of a fixed morphism $Y\lra X$ and then apply 
Proposition (\ref{bound3}) to bound the radical index. 
For this analysis we may assume that we are working over the
algebraic closure $\ov{K}$ of $K$.

Consider the map ${\widetilde \rho}: G\times \P(V) \lra \P(V)\times \P(V)$
defined by ${\widetilde \rho}=(\rho\,,\,{\rm pr}_2)$, where $\rho$ is
the action map and ${\rm pr}_2$ is the second projection. Let 
${\cal  Y}$  be the schematic image of ${\widetilde \rho}$. In this case 
${\cal  Y}$
gets the reduced induced scheme structure from the product and its
points are the closure of the image of the map ${\widetilde \rho}$.
We have the map $h: {\cal  Y} \lra \P(V)$ which is the composition of the
inclusion map to $\P(V)\times \P(V)$ with the second projection.
Let ${\cal Y}_x$ be the fiber of $h$ at a point $x \in \P(V)$. One observes
that  ${\cal Y}_x$ contains the closure of the orbit $O(x)$ of $x$. 

For any integral locally closed subscheme $Z \subset \P(V)$ we have the 
restriction of the
map ${\widetilde \rho}$ which we denote by ${\widetilde \rho}^Z$
from $G\times Z \lra \P(V)\times Z $.  Let ${\cal Y}^Z$ be its
schematic image. We again get the induced map $h^Z:{\cal Y}^Z \lra Z$
and we denote by ${\cal Y}^Z_x$ the fiber of the surjective map
$h^Z$ at $x\in Z$. 

Note that since $Z$ is reduced we have an open
subscheme $U_1 \subset Z$ where the map $h^Z$ is flat. One also
observes that the schematic image  of  
${\widetilde \rho}^Z |_{(G \times U_1)}$ is exactly 
$(h^Z)^{-1}(U_1)$.
Hence we can restrict the setup to $U_1$.

Since the actual image set ${\cal Z}'={\widetilde \rho}^Z(G \times U_1)$ 
is constructible, it contains an open subset 
$U' \subset {\cal Y}^{U_1}$. The map $U' \lra U_1$ is flat, hence if we
define $U=U'\,\cap \, U_1$ and restrict the whole setup to $U$, the 
image set ${\cal Z}={\widetilde \rho}^Z(G \times U)$
contains an open subset $U' \subset Y^{U}$ which maps surjectively onto
$U$. Now we can further translate the open set $U'$ by an element of
$g$ which acts on the first factor to obtain that ${\cal Z}$ is open
in ${\cal Y}^U$.
The upshot of this analysis is that for any subvariety $Z \subset \P(V)$  
  we can find an open subscheme $U \subset Z$ such that 
we can define an open reduced subscheme  ${\cal Z}^U$ of ${\cal Y}^U$ 
whose points are exactly the image of the map ${\widetilde \rho}^U$.
This also has the property that $O(x)$ is exactly the reduced part $({\cal Z}^U_x)_{\rm red}$
of the fiber ${\cal Z}^U_x$ at $x\in U$.

Let $0< l < {\rm dim}(V)$ be an integer. Let ${\rm Gr}_l(V)$ be the Grassmannian
of $l$-dimensional planes in $\P(V)$. We have a universal subscheme  
$i: H_l \subset {\rm Gr}_l(V) \times \P(V)$. This has the property that for
each $y\in {\rm Gr}_l(V)$ the inverse image 
$(H_l)_y=i^{-1}( \{y\}\times\P(V))\subset \P(V)$ is $\P(W_y)$, 
where $W_y$ is the $l$ dimensional 
subspace of $V$ associated to the point $y$.  

Now we consider the relative version of the setup to get a suitable         
intersection of the orbits with projective subspaces of $\P(V)$. Let $Z$ and $U$ be defined as before.

Let 
${\widetilde \rho}^U_l :{\rm Gr}_l(V) \times G \times U \lra 
{\rm Gr}_l(V) \times \P(V)\times U $ be
the map defined by $({\rm id}\,,\,{\widetilde \rho}^U)$. Inside 
${\rm Gr}_l(V) \times \P(V)\times U $ we have two subschemes namely 
${\rm Gr}_l(V) \times {\cal Z}^U$ and $ H_l \times U$. Let $Y$ be the
scheme theoretic intersection of these schemes. We have a map 
$f :Y \lra X={\rm Gr}_l(V) \times U $ which is obtained by composing the
inclusion map with the product of the first and the last projection map.

The map $f$ has the property that for any $x\in U$ and $y\in {\rm Gr}_l(V)$, 
$M_{y,x}'=f^{-1}(y\,,\,x)$ is the scheme theoretic intersection of
${\cal Z}^U_x$ and $\P(W_y)$. Hence by Proposition (\ref{bound3}), there 
exists an integer $n$ such that ${\rm Ri}(M_{y,x}') \le n$, for each 
$(y,\,x)\in X$. 

If $M_{y,x}$ is the schematic intersection of $O(x)$ and
$M_{y,x}'$ then we have an inclusion 
$(M_{y,x})_{{\rm red}} \subset M_{y,x} \subset M_{y,x}'$. Hence we
conclude that for each $(y\,,\, x) \in {\rm Gr}_l(V) \times U$, we have 
${\rm Ri}(\P(W_y) \cap O(x)) \le n$.

Take $Z \,= \,\P(V)$ and we obtain an open subscheme $U$ with a
bound $n$ for the radical index of $\P(W_y) \cap O(x)$ for every $x
\in U$ and $y\in {\rm Gr}_l(V)$. Then take the
complement of $U$ in $\P(V)$ and so on by induction there exists an
integer $n_l$ such that ${\rm Ri}(\P(W_y) \cap O(x)) \le n_l$ for each
$y\in {\rm Gr}_l(V)$ and $x\in \P(V)$.

We choose $N_2\,= \,{\rm Max}_{l=1}^r\{n_l\}$. Then for any $m \in M(K_s)$ which is 
not semistable and any choice of $x_m$ as before 
the subscheme $M(P)_{x_m}$ 
occurs as one such $M_{y,x}$ for some $y \in {\rm Gr}_l(V)$ for some
$l$ and $x=x_m$. This proves that for all choices of $x_m$ and $m$ 
we have ${\rm Ri}(M(P)_{x_m})\le N_2$.

As a next step in the proof of the Proposition (\ref{independence}) we
will show that for any $m\in M(K)$ the $K_s$-scheme  $M(P)_{x_m}$ 
admits a $K_s^{p^{-N_2}}$-rational point. For this we can first find an affine open subscheme
$U_{x_m}$ of $ M(P)_{x_m}$ which is defined over $K_s$. 
As the radical index of $U_{x_m}$ is $\le N_2$, the existence of rational points
on $U_{x_m}$ follows from the Lemma below.

\begin{lemma}\label{algebra1} 
Let $A$ be an affine $K_s$-algebra with radical index $\le p^n$. Then 
$A$ admits a $K_s^{p^{-n}}$-rational point.
\end{lemma}
\proof
We will denote by $A_n$ the $K_s^{p^{-n}}$-algebra $A\otimes _{K_s}\, K_s^{p^{-n}}$
and by ${\ov A}$ the ${\ov K}_s$-algebra defined by 
$A\otimes _{K_s}\, {\ov K}_s$. Since the radical 
${\rm Rad}(A)\otimes _{K_s}\, {\ov K}_s \subset {\rm Rad}({\ov A})$
we may assume that $A$ is reduced. We show that the natural inclusion
$$
{\rm Rad}(A_n)\otimes _{K_s}\,{\ov K}_s \subset {\rm Rad}({\ov A})
$$
is an isomorphism. This will prove that the $K_s^{p^{-n}}$-algebra 
$A_n/{\rm Rad}(A_n)$ is absolutely reduced and hence will admit a 
 $K_s^{p^{-n}}$-rational point. 

Let $f\in {\rm Rad}({\ov A})$. Then we can write 
$f= \sum _{i=1}^lf_i' \otimes a_i$ with $f_i \in A$ and $a_i \in {\ov
  K}_s$. 
If each of $a_i \in K_s^{p^{-n}}$ then already we have 
$f\in {\rm Rad}(A_n)\otimes _{K_s}\,{\ov K}_s$. 

Let $n_1\,>\,n$ be so chosen
that each of $a_i \in K_s^{p^{-{n_1}}}$. Using the identity  
${\ov A}= A_n \otimes _{K_s^{p^{-n}}}\, {\ov K}_s$, we have
an expansion of $f= \sum _{i=1}^j f_i \otimes b_i$ where $f_i \in A_n$ and 
$b_i \in K_s^{p^{-{n_1}}}$ such that $\{b_1,\ldots,b_j\}$ is linearly independent 
over  $K_s^{p^{-n}}$. 

We claim that $b_i^{p^n}$'s are linearly independent over $K_s$.
Suppose that $\sum c_i\,b_i^{p^n} =0$ with $c_i \in K_s$. 
Let $d_i \in K_s^{p^{-n}}$ be such that $d_i^{p^n}= c_i$. The above 
implies that $(\sum d_i\,b_i)^{p^n}\, =\,0$. This proves the claim.

The radical index of $A$ is $n$ we have $f^{p^n}=0$ 
and this gives us $0= f^{p^n}= \sum _{i=1}^j f_i^{p^n} \otimes b_i^{p^n}$.
Since $b_i^{p^n}$'s are linearly
independent over $K_s$ and $A$ is flat over 
$K_s$, we have
$f_i^{p^n}=0$. Hence $f_i \in {\rm Rad}(A_n)$. 
This proves the lemma (\ref{algebra1}).  $\hfill \square$

We have the orbit morphism $G \lra O(m)$ which is defined over
$K_s$. Also we have a $K_s^{p^{-N_2}}$-valued point of $M(P)_{x_m}$ which
is $K_s$-subscheme of $O(m)$. Hence we obtain a $K_s^{p^{-N_2}}$-valued
point $y$ of $O(m)$. Let $G_{m}$ be the isotropy subgroup scheme of
$G$ at $m$. Then we have a $K_s$ isomorphism $G/G_{m} \cong O(m)$.
Let $N_1$ be as defined in Proposition (\ref{bound}). 
We show that if $N=N_1\,+\,N_2$ then the instability parabolic $P(m)$ of
$m$ is defined over $K_s^{p^{-N}}$. For this we first show that there
is a $K_s^{p^{-N}}$-valued point $h$ of $G$ such that it maps to $y$. 
This statement follows from the fact that if $Y$
is the fiber of the map $G \lra O(m)$ at the point $y$, then $Y$
defines a principal $G_m$-bundle over $S={\rm Spec}(K_s^{p^{-N_2}})$.
Now it follows from Proposition (\ref{bound}) that the scheme  
$Y$ is an finite type affine scheme over $S$ whose 
radical index is $\le N_1$. Hence by lemma (\ref{algebra1}) we conclude
that $Y$ admits a $K_s^{p^{-N}}$-rational point. This proves that there is a 
 $K_s^{p^{-N}}$-valued point $h$ of $G$ such that $h\, m=y$.

The Lemma (\ref{parabolic}) now implies that 
the instability parabolic $P(m)$ of
$m$ is defined over $K_s^{p^{-N}}$ and this completes the proof of the Proposition 
(\ref{independence}). $\hfill \square$

\section{The proof of the main Theorems} 

In this Section we will prove the main Theorems (\ref{main}) and (\ref{sun})
stated in the introduction.
The basic strategy of the proof is to assume (\ref{main}) for the case of
lower semisimple rank groups and prove (\ref{sun}). Finally prove
the Theorem (\ref{main}) using (\ref{sun}) and the 
Proposition (\ref{independence}).

\begin{center}
{\bf Proof of Theorem (\ref{sun})} 

(Assuming Theorem (\ref{main}) for lower semisimple rank)
\end{center}

We fix the Borel subgroup $B$ and a maximal torus 
$T \subset B$ and the root datum as before.
Let $\Delta$ denote the set of simple roots.
Let 
$Q_{\alpha}$ be the maximal parabolic subgroup of $G$
containing $B$ corresponding to the simple root $\alpha \in
\Delta$.
We will denote by ${\mf q}_{\alpha}$ and ${\mf g}$ the lie algebras of 
$Q_{\alpha}$ and $G$ respectively.
Let $P$ be a parabolic subgroup of $G$ containing $B$. Let $L$ be its
Levi quotient.
We will use the following lemma which is proved in Biswas-Gomez 
(see proof of Theorem 4.1, in page 783, \cite{Biswas-Gomez}) . 

\begin{lemma}\label{biswas}
Let $P\subset Q_{\alpha}$ be an inclusion of parabolic
subgroups. There exists a filtration 
$$
0\,=\, V^{\alpha}_0\subset V^{\alpha}_1 \subset \ldots \subset 
V^{\alpha}_{a_{\alpha}}= {\mf g}/{\mf q}_{\alpha}
$$ 
of $P$ modules such that the following holds.
\begin{enumerate}
\item The unipotent radical $R_u(P)$ acts trivially on each successive 
quotients $W^{\alpha}_j=V^{\alpha}_j/V^{\alpha}_{j-1}$ and (hence) 
the identity connected component of the center of the Levi quotient $L$ acts by a scalar on the induced 
representation (denoted by $\rho ^{\alpha}_j$) of $L$.

\item The character $\chi ^{\alpha}_j$ of $P$ on the action of $P$ on
  $W_j^{\alpha}$ has the property that its
restriction to the maximal torus $T$ can be written as a non-positive 
linear combination of simple roots 
$\sum _{\beta \in \Delta} n^{\alpha}_{j,\beta} \beta$ 
with $n^{\alpha}_{j,\beta}\le 0$ and $n^{\alpha}_{j,\alpha}< 0$.
\end{enumerate}
\end{lemma}

Let $\Pi$ be the subset of simple roots defined by the property that 
$\alpha \in \Pi$ if and only if $P\subset Q_{\alpha}$. Here we can look at $P$ as $P_{\Pi}$.

For an $\alpha \in \Pi$, let $\chi^{\alpha}_0$ be the character of $P$ 
which is defined by the representation of $P$ on ${\mf g}/{\mf
q}_{\alpha}$. One can check that  the restriction of $\chi^{\alpha}_0$ to the maximal torus is a non-positive linear 
combination of simple roots with the coefficients of $\alpha$ being negative.
Let $\chi _P$ be
the character of $P$ defined by the representation of $P$ on 
${\mf g}/{\mf p}$. 

Let ${\mathcal T}$ be the finite set $\Pi _{\alpha \in \Pi}[0\,,\,a_{\alpha})$.
For an element ${\ov z} \in {\mathcal T}$,
the Lemma (\ref{biswas}) and the observation above imply that there
exists positive integers 
$n({\ov z})$ and $m^{\alpha}({\ov z})$ (for each $\alpha \in \Pi$) with
the property that the restriction of the  character 
$n({\ov z}) \chi _P \,-\,\sum_{\alpha \in \Pi}
m^{\alpha}({\ov z}) \chi ^{\alpha}_{{\ov z}_{\alpha}}$ to the maximal
torus $T$ is a linear combination of simple roots 
in $\Delta - \Pi$. 
This automatically implies that
\begin{equation}\label{char}
n({\ov z}) \chi _P \,=\,\sum_{\alpha \in \Pi}
m^{\alpha}({\ov z}) \chi ^{\alpha}_{{\ov z}_{\alpha}}
\end{equation}

We will define the constant $N_P$ by
\begin{equation}\label{numbers}
N_P\,=\, {\rm Max}_{{\ov z}\in T}
\{ \frac{1}{n({\ov z})} \sum _{\alpha \in \Pi} m^{\alpha}({\ov z})\}
\end{equation}

Lemma (\ref{biswas}) also implies that representation $\rho ^{\alpha}_j$
of $L$ takes the identity connected component of the center of $L$ to
connected component of the center of $GL(W_j^{\alpha})$. 
Hence by the induction assumption (on Theorem \ref{main}) 
there exists a constant
$C^L(X\,,\, \rho ^{\alpha}_j)$ such that for any rational semistable
$L$-bundle $E_L$ we have 
\begin{equation}\label{levi}
\mu_{{\rm max}}(E_L(W_j^{\alpha}))\,-\,
\mu_{{\rm min}}(E_L(W_j^{\alpha})) \le
C^L(X,\rho^{\alpha}_j) 
\end{equation}
Let $C_P={\rm  Max} \,\{C^L(X,\rho^{\alpha}_j)\} $, where the maximum is
taken over all $\alpha \in \Pi$ and $1\le j < a_{\alpha}$.

We also define the constant $M_P$ by setting 
$M_P ={\rm Max} _{\alpha \in \Pi,j}\{{\rm dim}(W^{\alpha}_j)\}$ 

Note that the constants $N_P$, $M_P$ and  $C_P$ depend on the parabolic 
$P$ and since there are only finitely many choices of parabolic
subgroups containing $B$, we will define  the constant 
$N$ ( respectively $M$ and $C$) to be the maximum of each $N_P$
(respectively $M_P$ and $C_P$) over all parabolics containing the
Borel subgroup $B$. 

We take a rational $G$-bundle $E$ over $X$. Let ${\rm Ideg}(E)$ be
its instability degree. Let $(P'\,,\sigma')$ be a Harder-Narasimhan
reduction.
One notes here that the
reduction $\sigma '$ satisfies the properties (1) and (2) stated in Section 2.

Let $F=F^*(E)$ be the Frobenius pull-back of $E$. Let $\sigma $ be its
Harder-Narasimhan reduction to a parabolic $P$ containing $B$. 
We will denote by $F_{\sigma}$ the $P$-bundle defined by
$\sigma$ and $T_{\sigma}$ the tangent bundle along the fibers of
$X$. We will denote by $F_{\sigma, L}$ the $L$-bundle obtained by
extension of $F_{\sigma}$ to $L$. 

We need to bound the slope of $T_{\sigma}$ in terms of
the slope of $T_{\sigma '}$.

For each $\alpha \in \Pi$, we have an inclusion $P\subset Q_{\alpha}$. This
gives rise to a reduction $\sigma _{\alpha}$ of $F$ to $Q_{\alpha}$
and we will denote by $F_{\sigma _{\alpha}}$ the $Q_{\alpha}$-bundle 
determined by this reduction.

By Lemma (\ref{biswas}) we have representations of $P$ in ${\mf g}/{\mf q}_{\alpha}$ and 
representations $\rho ^{\alpha}_j$ for $1 \le j <a_{\alpha}$.
These give rise to a vector bundle $T_{\sigma_{\alpha}}$ and a
filtration  $F_{\sigma}(V^{\alpha}_j)$ of  $T_{\sigma_{\alpha}}$ with
the property that successive quotients are isomorphic to 
$F_{\sigma, L}(W^{\alpha}_j)$.

The bundle $F$ (being the Frobenius pull-back) admits a $p$-connection $\nabla$ satisfying the properties
defined in Proposition (\ref{connection}).

We apply this to the reduction $\sigma _{\alpha}$ to get the map of
vector bundles 
$\nabla _{\sigma_{\alpha}}: T_X \lra T_{\sigma _{\alpha }}$.

First we consider the case when the above map is zero. Then there is a
reduction ${\widetilde \sigma}$ of $E$ to $Q_{\alpha}$ such that 
 $\sigma _{\alpha}=F^*({\widetilde \sigma})$. This has the effect
that ${\rm deg}(T_{\sigma _{\alpha}})= 
p \, {\rm deg}(T_{\widetilde \sigma})$.
Hence we have the inequality 
\begin{equation}\label{0th}
{\rm deg}(T_{\sigma _{\alpha}})\,\ge \, p \, {\rm Ideg}(E)
\end{equation}

Now suppose that map $\nabla _{\sigma_{\alpha}}$ is not zero. Then
there is a $j$ such that the image of $\nabla _{\sigma_{\alpha}}$ is contained in $F_{\sigma}(V^{\alpha}_j)$
 and not in $F_{\sigma}(V^{\alpha}_{j-1})$.
Hence we get a non-trivial map 
$T_X \lra F_{\sigma, L}(W_j^{\alpha})$.
This implies that $\mu_{{\rm min}}(T_X)\le 
\mu_{\rm max}(F_{\sigma, L}(W_j^{\alpha}))$.

Let $C_X =\mu_{{\rm min}}(T_X)$. This combined with equation
(\ref{levi}) and the fact that $C_X -C$ can be made to be negative 
implies that 
\begin{equation}\label{jth}
{\rm deg}(F_{\sigma, L}(W_j^{\alpha})) \ge (C_X -C)\, M
\end{equation}

Further the right hand side in the inequalities (\ref{0th}) and
(\ref{jth}), being negative, can be summed up to get a common 
right hand side, namely $ p\, {\rm Ideg}(E)\,+\,(C_X -C)\, M$

Hence for any $\alpha \in \Pi$ either the inequality (\ref{0th}) holds
or the inequality  (\ref{jth}) holds for some choice of $j$. This implies
that if we vary $\alpha \in \Pi$ we obtain an element 
${\ov z}\in {\mathcal T}$.
Hence for the element ${\ov z}$ using the formula (\ref{char})
we get 
$$ n({\ov z}){\rm deg}(T_{\sigma})\,\geq\,
\sum_{\alpha \in \Pi}
(m^{\alpha}({\ov z}))(p \,{\rm Ideg}(E)\,+\,(C_X -C)\,M)
$$ 
This implies that ${\rm deg}(T_{\sigma}) 
\ge N(p\,{\rm Ideg}(E)\,+\,(C_X -C)\, M)$. This completes the proof
of Theorem (\ref{sun}). $\hfill \square$

The above theorem and an induction argument also proves the following.

\begin{corollary}\label{remark} 
There exists  
constants $C$ and $N$ (independent of $E$) such that 
$${\rm Ideg}((F^n)^*E) \ge  p^n\, N\,{\rm Ideg}(E)\,+\,C$$.
\end{corollary}

\begin{center}
{\bf Proof of the Theorem (\ref{main})}
\end{center}

We fix a Borel subgroup $B_1$ of ${\rm GL}(V)$. For a parabolic 
$P_1$ of  ${\rm GL}(V)$ containing $B_1$, we have an action of $G$ on 
$M_{P_1}={\rm GL}(V)/P_1$. We fix a
representation  ${\rm GL}(V) \lra {\rm GL}(V_{P_1})$ such that it
defines an embedding of ${\rm GL}(V)/P_1 \subset \P(V_{P_1})$ with the
property that the character of $P_1$ on $V_{P_1}$ is a positive 
multiple $m_{P_1}$ of the character of $\chi_{P_1}$ associated to the
restriction of the adjoint representation of $P_1$ on the vector space 
${\mf gl}(V)/{\mf p}_1$.

The line bundle ${\cal O}(-1)$ on $\P(V_{P_1})$ when restricted to 
${\rm GL}(V)/P_1$ defines an anti ample line bundle ${\cal L}_{P_1}^{-1}$ which is
also defined by the character $-m_{P_1}\,\chi_{P_1}$. 

For a rational $G$-bundle $E$ over $X$ we have a rational fiber bundle
$E(M_{P_1})$ and the line bundle ${\cal O}(-1)$ gives a
rational line bundle $E({\cal L}_{P_1}^{-1})$. 

Let $x_0$ be the generic point of $X$. Let $E(G)_0$ the group scheme 
over $K=k(X)$ associated to $E$ at 
the generic point $x_0$. Then we have an action of 
$E(G)_0$ on $E(M_{P_1})_0$ which is linearized with respect to the 
line bundle $E({\cal L}_{P_1})_0$ over $K$.

\begin{lemma}\label{indipendence}There exists a constant $N$,
depending only on $G$ and $X$, such that for any rational $G$-bundle $E$ and for any parabolic $P_1$ containing $B_1$ the  
instability parabolic for any $K$-valued non-semistable point of $E(M_{P_1})$ (for
the above action of $E(G)_0$) is defined over $K^{p^{-N}}$.
\end{lemma}
\proof

Let $E_0$ be the principal $G$-bundle over $K$ obtained by restriction
of $E$ to the generic point of $X$. One observes that $E_0$ becomes 
trivial over a finite separable
extension of $K$, hence when we change the base to $K_s$,
the separable closure of $K$, we get an isomorphism 
$E_0 \otimes _KK_s \cong G \otimes _k K_s$. This isomorphism now
canonically extends to give an isomorphism of 
$E(G)_0 \otimes _KK_s \cong G \otimes _k K_s $ and 
$E(M_{P_1})_0 \otimes _KK_s \cong M_{P_1}\otimes _k K_s$, and 
the last one being compatible with group actions, 
and also of the isomorphisms between the ample line bundles on these spaces.

For the induced action of $G \otimes _k K_s$ on  $M_{P_1}\otimes _k K_s$
which is linearised by ${\cal L}_{P_1}\otimes _k K_s$, by 
Proposition (\ref{independence}), it follows that there is a positive
integer $N$ such that for any non-semistable point $m$ of 
$M_{P_1}\otimes _k K_s$ the instability is defined over
$K_s^{p^{-N}}$. Since there are only finitely many parabolic subgroups containing $B_1$ we can find a constant $N$
which works for all these parabolic subgroups.

This implies that the group scheme $E(G)_0 \otimes _KK_s$ for any $E$ also has
the same property. The Galois descent argument implies that 
instability parabolic of a $K$-valued point of $E(M_{P_1})_0$ is defined
over the field extension $K^{p^{-N}}$, with $N$ being independent of
the rational $G$-bundle $E$ and the reduction ${\sigma}$. 
This proves the lemma.$\hfill \square$

A rational reduction $\sigma $ of $\rho_*E$ to $P_1$ gives a $K$-rational 
point $\sigma (x_0)$ of $E(M_{P_1})_0$.
If this point is semistable then by Proposition (\ref{RR1}), we have
\begin{equation}\label{funny}
{\rm deg}(\sigma^*E({\cal L}_{P_1}))\,= \, m_{P_1}{\rm deg}(T_{\sigma})
\geq 0
\end{equation}

If the point $\sigma(x_0)$ is not semistable we have an integer $N$ 
prescribed by Lemma (\ref{indipendence})
such that its instability parabolic $P'_0$ is defined over $K^{p^{-N}}$.

One observes that pull-back by the $N$-th Frobenius
morphism $F^N$ of $X$, 
the action of the generic fibre $((F^N)^*E(G))_0=(F^N)^*_K(E(G)_0)$ on 
$((F^N)^*E(M_{P_1}))_0=(F^N)^*_K(E(M)_0)$ is the base
change by the Frobenius $F^N_K:Spec(K)\to Spec(K)$ of the action of
$E(G)_0$ on $E(M)_0$.  Hence the point $((F^N)^*\sigma)(x_0)$ 
has an
instability parabolic $P'_0$ for the action of $(F^N)_K^*(E(G)_0)$ defined over
$K^{p^{-N}}$ (see proof of Theorem (3.23) of \cite{Ramanan-Ramanathan}).

The parabolic subgroup $P'_0$ defines a rational reduction of the structure 
group $\tau$ of $(F^N)^*E$ to a
parabolic subgroup $P'\subset G$ with $P'_0=(E_{\tau}(P'))_0$.  

Since $P'_0$ is defined over $K^{p^{-N}}$, by Remark (\ref{uniqueness}), 
the instability 1-PS for $\sigma (x_0)$ is also defined over $K^{p^{-N}}$.

The Proposition (\ref{RR2}) (applied to bundle $(F^N)^*E$ and the point  $((F^N)^*\sigma)(x_0)$) implies that there is a positive integer 
$r$ and a
dominant character $\chi$ of $P'$ such that the following inequality
holds.
\begin{equation}\label{eqn}
-r \, {\rm deg}({((F^N)^*\sigma)}^*(((F^N)^*E)({\cal L}_{P_1}))) \,\le\, {\rm deg}(\chi _* (((F^N)^*E)_{\tau})) 
\end{equation}  

Let $\Pi \subset \Delta$ be the subset defining the parabolic $P'$.
Let $Q_{\alpha}$ be the maximal parabolic subgroup of $G$ containing
$P'$ defined by $\alpha$. 
Let $\chi _{\alpha}$ be the dominant character of $Q_{\alpha}$ 
defined by the representation ${\mf g}/{\mf q}_{\alpha}$.
There is a positive integer $m_{\alpha}$ such that  
we have $\chi _{\alpha}|_T=-m_{\alpha}w_{\alpha}$ where  
$w_{\alpha}$ are the fundamental weights of $G$ with respect to a fixed 
maximal torus contained in $P'$.

Let $L' = P'/R_u(P')$ and $Z_0(L')$ be the connected component of the
center of the Levi $L'$.

By Lemma (\ref{git4}) we have a character 
$\chi '$ of the maximal torus $T \subset P'$ such that 
$\chi'|_{Z^0(L')}=\chi |_{Z^0(L')}$
and $\chi' =r\,\nu(v\,,\lambda)\,||\lambda || \,l_{\lambda}$. Here
$l_{\lambda}$ is the dual of $\lambda$. 

The above equality can be rewritten as 
$\chi' = r\, \,\nu(\sigma_N\,,\, \lambda) l_a$ where 
$a \in {\cal X}_*(T)\otimes \Q$ is the
element in the unit sphere defined by $a=\lambda/\|\lambda\|$.

By Lemma (\ref{git1}) 
there exists a constant $B_G$ such that 
for every point $m \in M_{P_1}\otimes _k {\ov K}$ and 
$\lambda \in {\cal X_*(T)}\otimes \Q$ we have 
\begin{equation}\label{git}
\nu(m\,,\,\lambda) \,\le \, B_G
\end{equation}

The Weyl group invariant scalar product 
on $ {\cal X}_*(T)\otimes \Q$ induces a scalar product on 
$ {\cal X}^*(T)\otimes \Q$.
The following lemma is an elementary calculation. 
\begin{lemma}\label{basic}
Let $S = \{l \in {\cal X}^*(T)\otimes \Q \,\,|\,\,\|l\|=1 \}$.
Then there exists a constant $A_G$ such that
for each $l\in S$ if $l=\sum_{\alpha \in \Delta}r_{\alpha}w_{\alpha}$
then $|r_{\alpha}|\,\le \,A_G$ for each $\,\alpha \in \Delta$.
\end{lemma}

One notes that under the scalar product we have 
$\|l_{a}\| =1$. Since $l_{a}$ is trivial on the center of
$G$, the Lemma (\ref{basic}) implies that 
$l_a= \sum _{\alpha \in \Delta}r_{\alpha}w_{\alpha}$ with 
$|r_{\alpha}| \,\le\, A_G$.

This along with the above description of $r$ and $\chi$ we get 
\begin{equation}
\chi |_{Z_0(M')} \,=
\, r \,\nu(\sigma_N\,,\, \lambda) \, 
\sum _{\alpha \in \Pi}r_{\alpha}w_{\alpha}
\end{equation}
the last equality can be rewritten in terms of $\chi _{\alpha}$ as follows.
$$-\chi \,=\,  r\,\nu(\sigma_N\,,\, \lambda) \, 
\sum _{\alpha \in \Pi}(r_{\alpha}/m_{\alpha})\chi_{\alpha}
$$  

Using the fact that ${\rm deg}({\chi_{\alpha}}_{*}(((F^N)^*E)_{\tau}))\,\ge \, {\rm
 Ideg}((F^N)^*E)$
and combining it with the inequality (\ref{git}) 
we obtain 
$$
{\rm deg}(\chi _* (((F^N)^*E)_{\tau}))\le  -r\,B_G \,A_G \,|\Delta|  {\rm Ideg}((F^N)^*E).
$$

This along with (\ref{eqn}) implies that there exists a constant  $C(G)$ depending only on $G$
such that 
${\rm deg}(((F^N)^*(\sigma))^*E({\cal L}_{P_1})) \,\ge C(G) {\rm Ideg}((F^N)^*E)$. 

Since $((F^N)^*{\sigma})^*(((F^N)^*E)({\cal L}_{P_1}))= (\sigma ^*(E({\cal L}_{P_1})))^{p^N}$ 
it follows that 
${\rm deg}(\sigma^*(E({\cal L}_{P_1})) \,\ge (C(G)/p^N) {\rm Ideg}((F^N)^*E)$. By
(\ref{funny}) we have  ${\rm deg}T_{\sigma} \ge (C(G)/(m_{P_1}\,p^N)){\rm Ideg}((F^N)^*E)$
for every rational reduction $\sigma$ of $\rho_*E$ to $P_1$ which has the property that $\sigma(x_0)$ is not semistable.

This implies that we have a constant $C(G\,, \rho)$ such that 
$$
{\rm Ideg}(\rho_*E) \ge C(G\,,\rho) \frac{{\rm Ideg}((F^N)^*E)}{p^N}
$$
By Corollary (\ref{remark}) and Lemma (\ref{GLN}) we are through with the proof of the
Theorem (\ref{main}). $\hfill \square$

\rem \label{general}  The proof actually shows the following more general result.
Let $\rho: G\to G'$ be a homomorphism of connected reductive groups which takes the 
identity connected component of the center of $G$ to the center of $G'$. Then 
there exists constants $C$ and $C'$ (depending only on $X$, $G$, and $\rho$)
such that for any rational $G$-bundle $E$ over $X$ we have 
$$
{\rm Ideg}(\rho_*E) \ge C \,{\rm Ideg}(E) \,+\,C'.
$$
With this formulation, Theorem (\ref{sun}) is a special case of this result when applied to
the Frobenius homomorphism of $G$.

\section{Boundedness of semistable bundles} 
In this section we prove the boundedness of semistable $G$-bundles on
$X$ under the assumption stated in the introduction. 
From now on we work with $G$-bundles on $X$ and not the
rational ones.

Let ${\cal X}(G)$ be the group of characters of $G$.
Let $A^k(X)$ be the $k$-th Chow group. To a principal $G$-bundle $E$,
recall the definition of the degree  
$d_E\in {\rm Hom}({\cal  X}(G)\,,\,A^1(X))$ of a principal $G$ bundle
$E$ from the introduction.

We fix a collection of elements $c_i \in A^i(X)$ for $2 \,\le\, i \,\le \,n={\rm dim}(X)$ 
and also we fix an element 
$d \in {\rm Hom}({\cal  X}(G)\,,\,A^1(X))$. 

Under the assumptions in 
Theorem (\ref{boundedness}), we will show that  
the set ${\cal S}_G(d\,;\,c_2 \,,\ldots, c_{n})$ 
of isomorphism classes of semistable $G$-bundles 
$\{E\}$ with degree $d_E=d$ and the Chern classes 
$c_i({\rm ad}(E))=c_i$ is bounded. Our proof is based on Proposition (4.12)
of \cite{Ramanathan}

We begin with an elementary lemma which allows us to use
representations.
\begin{lemma}\label{rep}
There is a faithful completely reducible rational representation of $G$.
\end{lemma}
\proof
For any irreducible representation $\rho$ of $G$ in a vector space $V$ 
let ${\rm ker}({\rho})$ be the kernel. 
We first show that 
$N=\bigcap _{\rho}{\rm ker}({\rho})$, over all irreducible
representations, is
trivial. This is because if $\rho_1$ is a faithful representation of
$G$ (hence of $N$) on a vector space $W$ then there is a filtration of 
$W$ such that successive quotients are irreducible. 
This implies that  $\rho_1(N) \subset {\rm GL}(W)$ lies in a
parabolic subgroup and its image when composed with the projection to
the Levi quotient is trivial. Hence $N$ is a unipotent
normal subgroup scheme. Let $N_0$ be the identity component of $N$.
Since $G$ is connected, using the conjugation map $G \times N_0 \lra N$
defined by $(g\,,n)\mapsto gng^{-1}$, we check that $N_0$ is normal.
Again using the conjugation map, this time from 
$G \times (N_0)_{\rm red} \lra N$ we see that $(N_0)_{\rm red}$ is
also normal. Since $G$ is reductive this proves that $N$ is a finite
subgroup scheme of $G$. Now using the conjugation map for the third
time we get that $N$ is central and hence it is diagonalizable. Now we see that the
representation $\rho _1$ restricted to $N$ is trivial which is a 
contradiction.

Now by dimension and length count we can find finitely many 
irreducible representations $\rho _i$, for $i= 1\,,\ldots ,m$, of $G$ 
such that $\bigcap _{i=1}^m{\rm ker}(\rho _i) =0$. This proves the lemma. $\hfill \square$

We also have another general lemma.

\begin{lemma}\label{rep2}
Let $\rho :{\rm GL}(V) \lra {\rm GL}(W)$ be a representation of 
${\rm GL}(V)$. Let $E_1$ and $E_2$ be two ${\rm GL}(V)$ bundles over
$X$ such that $c_i(E_1)=c_i(E_2)$ for each $i$. Then we have
$c_i(E_1(W))=c_i(E_2(W))$ for each $i$. 
\end{lemma}
\proof 
Let  $A^*({\rm BGL}(V))$ (respectively $A^*({\rm BGL}(W))$) be the Chow 
ring of $ {\rm BGL}(V)$ (respectively $ {\rm BGL}(W)$). Then one knows that 
$A^*({\rm BGL}(V))\cong \Z[c_1\,,\ldots ,c_n]$(see \cite{Totaro}). 
The representation $\rho$ gives rise to the map 
$\rho ^*:A^*({\rm BGL}(W))\lra A^*({\rm BGL}(V))$.   
Now we have the classifying maps $f_{E_i}: X \lra {\rm BGL}(V)$ for $i=1\,,2$.
The conditions of the lemma imply that the induced maps 
$f_{E_i}^* :A^*({\rm BGL}(V))\lra A^*(X)$ are equal. This implies that 
the maps $f_{E_1}^*\circ \rho ^* = f_{E_2}^*\circ \rho ^*$. Hence
the Lemma follows. $\hfill \square$

We continue with the proof of the Theorem (\ref{boundedness}).
By lemma (\ref{rep}), we have a completely reducible representation 
$\rho =\bigoplus \rho_i $ on $V = \bigoplus V_i$ of $G$ which is
faithful. Here $V_i$ are the irreducible components.
 
For a fixed $i$, the representation $\rho_i$ takes the identity connected component
of the center of $G$ to the center of ${\rm GL}(V_i)$. 
Hence for $E$ and $E'$ in ${\cal S}_G(d\,;\,c_2 \,,\ldots, c_{a})$, 
we have $c_1(E(V_i))=c_1(E'(V_i))$.

The representation $\rho_i$ induces 
a Lie algebra homomorphism ${\mf g}\lra {\rm End}(V_i)$. 
The lemma (\ref{rep2}) now proves that the two vector bundles 
$E({\rm End}(V_i)) ={\rm End}(E(V_i))$ and 
$E'({\rm End}(V_i)) ={\rm End}(E'(V_i))$
have same Chern classes.
Since the Chern classes of a vector bundle
are completely determined by the first Chern class and the 
Chern classes of the endomorphism bundle we conclude that 
the Chern classes of $E(V_i)$ are are independent of the 
individual members $E$ in ${\cal S}_G(d\,;\,c_2 \,,\ldots, c_{n})$.

By Theorem (\ref{main}) there are constants $C_i=C(X,\rho_i)$ for the
representation $\rho _i$. 

Now the assumptions in Theorem (\ref{boundedness}) imply that 
there is a  finite type scheme $S_i$ and a
family $U_i$ of vector bundles over $S_i \times X$ which contains every
member of 
${\cal S}_{C_i}(r_i\,;\,c_1(E(V_i))\,,\ldots c_{n}(E(V_i)))$ occurs. 

Now we take the scheme $S= \Pi_{i=1}^mS_i$
 and $U=U_1 \times_X U_2 \times_X \,\ldots \, U_m$. Then Theorem 
(\ref{boundedness}) follows from the arguments in the last part of the
 Proposition 3.1 in \cite{Holla-Narasimhan}. $\hfill \square$

\rem It should be possible to prove a version of the above
 Theorem for the case of principal $G$-sheaves in the sense of 
\cite{Gomez-Sols}.


\bigskip

{\it 
Address : \\
School of Mathematics, Tata Institute of Fundamental
Research, Homi Bhabha Road, Mumbai 400 005, India.
e-mail: yogi@math.tifr.res.in,\\
\\
S.I.S.S.A. Via Beirut 4, 34013 Trieste, Italy\\
e-mail: coiai@sissa.it 
} 
\bigskip

\end{document}